\documentclass{amsart}
\usepackage{amssymb, latexsym}
\newcommand{\half}{\frac{1}{2}}

\newcommand{\e}{\epsilon}

\newcommand{\lam}{\lambda}

\newcommand{\R}{{\mathbb R}}

\newcommand{\Lp}[2]{{\left\| #2 \right\|}_{L^ #1 }}

\newcommand{\Hsup}[2]{{{\left\| #2 \right\|}_{{H^#1}}}}

\newcommand{\X}[3]{{{\left\| #3 \right\|}_{{X_{#1,#2}}}}}
\newcommand{\Xd}[3]{{{\left\| #3 \right\|}_{{X^\delta_{#1,#2}}}}}

\newcommand{\Xsb}{{X_{s,b}}}

\theoremstyle{definition}

\theoremstyle{remark}
\newtheorem{remark}{Remark}
\theoremstyle{proposition}
\newtheorem{proposition}{Proposition}
\theoremstyle{lemma}
\newtheorem{lemma}{Lemma}
\theoremstyle{corollary}

\numberwithin{equation}{section}

\begin{document}

\title{Global well-posedness for KdV in Sobolev spaces of negative index}
\author{J.~Colliander}
\thanks{J.E.C. is supported in part by an N.S.F. Postdoctoral 
           Research Fellowship.}
\address{University of California, Berkeley}
\author{M.~Keel}
\thanks{M.K. is supported in part by N.S.F. Grant DMS 
                         9801558}
\address{Caltech}
\author{G.~Staffilani}
\thanks{G.S. is supported in part by N.S.F. Grant DMS 9800879 and by a
Terman Award.}
\address{Stanford University}
\author{H.~Takaoka}
\address{Hokkaido University}
\author{T.~Tao}
\thanks{T.T. is a Clay Prize Fellow and is supported in part by grants
from the Packard and Sloan Foundations.}
\address{University of California, Los Angeles}


\begin{abstract}

The initial value problem for the Korteweg-deVries equation on the line
is shown to be globally well-posed 
for rough data. In particular, we show global well-posedness for 
initial data in $H^s (\R )$ for $-\frac{3}{10} < s$.

\end{abstract}

\maketitle

\section{Introduction}

Consider the initial value problem for the Korteweg-deVries (KdV) equation
\begin{equation}
\left\{ \begin{matrix}
\partial_t u + \partial_x^3 u + \half \partial_x {(u^2)} = 0, &x \in \R,  \\
u(0) = \phi ,
\end{matrix}
\right.
\label{kdvivp}
\end{equation}
for rough initial data $\phi \in H^s ( \R) , ~s < 0$. This problem is known 
\cite{KPVBilin} to be locally well-posed provided $- \frac{3}{4} < s $. 
For $s \geq 0$,
the local result and $L^2$ norm conservation imply \eqref{kdvivp} is globally 
well-posed \cite{B1}. Recently, a direct adaptation \cite{CST99} of 
Bourgain's high-low frequency technique \cite{BAMS}, \cite{Brefine} 
showed \eqref{kdvivp} is
globally well-posed for $\phi \in H^s \cap {\dot{H}}^a $ for certain $s,a <0$.
A modification of the high-low frequency technique, first used in 
\cite{KeeTao98}, is presented in this paper which establishes global 
well-posedness of \eqref{kdvivp} in $H^s ( \R ), ~ -\frac{3}{10} < s $.

A subsequent paper \cite{CKSTTKdV} will establish 
that \eqref{kdvivp} is globally well-posed in $H^s ( \R )$
for $- \frac{3}{4} < s.$ The simplicity of the argument presented here
may extend more easily to other situations, such as in our treatment 
\cite{CKSTTNLS} of cubic $NLS$ on $\R^2$ and $NLS$ with derivative
in $\R$ \cite{CKSTTDNLS}.

\noindent{\bf{The Multiplier operator $I$}}

Let $s< 0$ and $N \gg 1$ be fixed. Define the Fourier multiplier operator
\begin{equation}
\label{Idefined}
{{\widehat{Iu}}}( \xi ) = m( \xi ) {\widehat{u}} ( \xi ),
~m( \xi ) = \left\{ \begin{matrix}
1, & |\xi | < N, \\
N^{-s} {{|\xi |}^s} , & |\xi | \geq 10 N
\end{matrix}
\right.
\end{equation}
with $m$ smooth and monotone.
The operator $I$ (barely) maps $H^s ( \R ) \longmapsto L^2 ( \R )$. Observe
that on low frequencies $\{ \xi : |\xi | < N \},~I$ is the identity operator.
Note also that $I$ commutes with differential operators. 
The operator $I^{-1}$ is the
Fourier multiplier operator with multiplier $\frac{1}{m(\xi )}$.

\noindent{\bf{An almost $L^2$ conservation property of \eqref{kdvivp}}}

Let $\phi \in H^s (\R ), ~ - \frac{3}{4} < s < 0$ in \eqref{kdvivp}. There
is a $\delta = \delta ( \Hsup s {\phi } ) > 0$ such that \eqref{kdvivp}
is well-posed for $t \in [0, \delta]$. We observe using the Fundamental 
Theorem of Calculus, the equation, and integration by parts that
\begin{eqnarray*}
{{\| Iu ( \delta ) \|}_{L^2}^2 }
&=& {{\| Iu (0 ) \|}_{L^2}^2 } + \int_0^\delta \frac{d}{d\tau }
( Iu( \tau) , Iu ( \tau ) ) d\tau, \\
& = &  {{\| Iu (0 ) \|}_{L^2}^2 } + 2 \int_0^\delta
( I {\dot {u}} ( \tau ) , Iu( \tau )) d \tau , \\
&= & {{\| Iu (0 ) \|}_{L^2}^2 } + 2 \int_0^\delta
( I ( - {u_{xxx}} - \half \partial_x [u^2] )(\tau) , Iu( \tau )) d\tau \\
& = &  {{\| Iu (0 ) \|}_{L^2}^2 } + \int_0^\delta
(I ( - \partial_x [u^2] ), Iu) d\tau.
\end{eqnarray*}
Finally, we add $0 = \int_0^\delta \int \partial_x ( {{I(u)}^2} ) I(u) d \tau$
to observe
\begin{equation}
\label{encounter}
{{\| Iu ( \delta ) \|}_{L^2}^2 }
= {{\| Iu (0 ) \|}_{L^2}^2 } + \int_0^\delta \int
\partial_x \left\{
{{(I(u))}^2} - I( u^2 ) \right\}~ Iu ~dx  d\tau.
\end{equation}
This last step enables us to take advantage of some internal cancellation.
We apply Cauchy-Schwarz as in \cite{Staff97} and bound the integral
above by
\begin{equation}
\Xd 0 {- \half -} { \partial_x \{ {{(I(u))}^2} - I( u^2 ) \} }
\Xd 0 {\half + } { Iu } .
\label{postcs}
\end{equation}

\begin{remark}
An effort to find a term providing more cancellation than
$\int_0^\delta \int \partial_x ( {{I(u)}^2} ) I(u) d \tau$ used above
led to the general procedure described in \cite{CKSTTKdV}.

\end{remark}

\begin{proposition}
(A variant of local well-posedness) The initial value problem
\eqref{kdvivp} is locally well-posed in the Banach space
$I^{-1} L^2 = \{ \phi \in H^s ~{\mbox{with norm}}~ \Lp 2 {I \phi } \}$.
with existence lifetime $\delta$ satisfying
\begin{equation}
\label{variantlifetime}
\delta \gtrsim {{ \| I \phi \|}_{L^2}^{-\alpha}}, ~{\mbox{for some}}~ 
\alpha > 0,
\end{equation}
and moreover
\begin{equation}
\label{variantspacetimenormbound}
\Xd 0 {\half +} {Iu} \leq C \Lp 2 {I \phi }.
\end{equation}
\end{proposition}

This proposition is not difficult to prove using the argument in
\cite{KPVBilin}. Using Duhamel's formula and $\Xsb$ space properties
reduces matters to proving the bilinear estimate
\begin{equation}
\X 0 {-\half+} {\partial_x I(u v ) } \leq C \X 0 {\half+ } {Iu }
\X 0 {\half+} {Iv}
\label{lwpbilvar}
\end{equation}
to obtain the contraction. The space-time norm bound is then implied
by the contraction estimate. The estimate \eqref{lwpbilvar} follows from the 
next proposition and the bilinear estimate of Kenig, Ponce and
Vega \cite{KPVBilin}.

\begin{proposition} (Extra smoothing)
The bilinear estimate
\begin{equation}
\label{smoothingest}
\Xd 0 {-\half -} { \partial_x \{ I(u) I(v) - I( uv ) \} }
\leq C N^{-\frac{3}{4}+} \Xd 0 {\half+ } {Iu }
\Xd 0 {\half + } {Iv }.
\end{equation}
holds.
\end{proposition}
Recall the bilinear estimate
$\X 0 {-\half + } {\partial_x (u v) } \leq C \X 0 {\half+ } u
\X 0 {\half+ } v $ from \cite{KPVBilin}. 
Proposition 2 reveals a smoothing beyond
the recovery of the first derivative for the particular quadratic
expression encountered above in \eqref{encounter}. We prove 
Proposition 2 in the next section.

The required pieces are now in place for us to give the proof of 
global well-posedness of \eqref{kdvivp} in $H^s ( \R) , ~ - \frac{3}{10} < s.$
Global well-posedness of \eqref{kdvivp} will
follow if we show well-posedness on $[0,T]$ for arbitrary $T>0$. We 
renormalize things a bit via scaling. If $u$ solves \eqref{kdvivp} then
$u_\lam ( x ,t ) = {{(\frac{1}{\lam}  )}^2} 
u( \frac{x}{\lam} , \frac{t}{\lam^3})  $ solves \eqref{kdvivp} with initial
data $\phi_\lam (x,t) = 
 {{(\frac{1}{\lam}  )}^2} 
\phi ( \frac{x}{\lam} ).$
Note that $u$ exists on $[0,T]$ if and only if $u_\lam$ exists on $[0,
 \lam^3 T]$. A calculation shows that
\begin{equation}
\Lp 2 {I \phi_\lam } \leq C {\lam^{-\frac{3}{2} - s }} N^{-s}  \Hsup s {\phi} .
\label{scaling}
\end{equation}
Here $N= N(T)$ will be selected later but we choose $\lam = 
\lam(N)$ right now by requiring
\begin{equation}
 C {\lam^{-\frac{3}{2} - s }} N^{-s}  \Hsup s {\phi} \thicksim 1 \implies
\lam \thicksim N^{- \frac{2s}{3 + 2s }}.
\label{lamchosen}
\end{equation}

We now drop the $\lam$ subscript on $\phi$ by assuming that 
\begin{equation}
\Lp 2 {I \phi} = \epsilon_0 \ll 1
\label{normalized}
\end{equation}
and our goal is to construct the solution of \eqref{kdvivp} on the time
interval $[0, \lam^3 T]$.

The local well-posedness result of Proposition 1 shows we can construct the
solution for $t \in [0, 1]$ if we choose $\epsilon_0$ small enough. 
The almost $L^2$ conservation property shows
${{\| I u(1) \|}_2^2} \leq {{\| Iu(0) \|}_2^2} + N^{-\frac{3}{4}+}
{{\| Iu \|}_{ {{X_{0, \half+}}} }^3}.$ Using \eqref{variantspacetimenormbound}
and \eqref{normalized} gives
\begin{equation*}
{{\| Iu(1) \|}_2^2} \leq \e_0^2 + N^{-\frac{3}{4}+}.
\end{equation*}
We can iterate this process $N^{\frac{3}{4}-}$ times before doubling 
$\Lp 2 {Iu(t)}$. Therefore, we advance the solution by taking 
$N^{\frac{3}{4}-}$ time steps of size $O(1)$. We now restrict 
$s$ by demanding that
\begin{equation}
N^{\frac{3}{4}-} \gtrsim \lam^3 T = N^{\frac{-6s}{3 + 2 s}} T
\label{scondition}
\end{equation}
is ensured for large enough $N$,
so $s> - \frac{3}{10}$.

\section{Proof of the bilinear smoothing estimate}

This section establishes Proposition 2. We distinguish the {\bf{very low
frequencies}} $\{ \xi : |\xi | \lesssim 1 \}$, the {\bf{low frequencies}}
$\{ \xi : 1 \lesssim |\xi | \lesssim \half N \}$ 
and the {\bf{high frequencies}}
$\{ \xi: \half N \lesssim |\xi | \}$. 
Decompose the factor $u$ in the bilinear estimate
by writing $u = u_{vl} + u_l + u_h $ with ${\widehat{u_l}}$ supported
on the low frequencies and similarly for the very low and high
frequency pieces. We decompose $v$ the same way. Since $I$ is the identity
operator on the low and very low frequencies, we can assume one of the 
factors $u,v$ in the estimate to be shown has its Fourier transform supported
in the high frequencies. Symmetry allows us to assume $u = u_h$ and
we need to consider the three possible interactions of $u_h$ 
with $v_{vl},~ v_l$ and $v_h$. Finally, since we are considering (weighted) 
$L^2$ norms, we can replace ${\widehat{u}}$ and ${\widehat{v}}$ by
$|{\widehat{u}}|$ and $|{\widehat{v}}|$. Assume therefore that $
{\widehat{u}}, {\widehat{v}} \geq 0$.

\noindent{\bf{Very low/high interaction}}

An explicit calculation shows that
\begin{equation}
{\mathcal{F}} \left( 
\partial_x \{ I( u_h v_{vl} ) - I(u_h) v_{vl} \} \right) (\xi)
= \int_{\xi = \xi_1 + \xi_2} i \xi [m(\xi ) - m( \xi_1 ) ]
{\widehat{u_h}}( \xi_1 ) {\widehat{v_{vl}}} ( \xi_2 ) ,
\label{expvlhi}
\end{equation}
where $\mathcal{F}$ denotes the Fourier transform.
The mean value theorem gives
\begin{equation*}
| m( \xi) - m( \xi_1 ) | \leq | m' ( {\tilde{ \xi_1 }} ) | |\xi_2 |,
\end{equation*}
which may be interpolated with the trivial estimate to give
\begin{equation}
\label{derivshuttle}
| m( \xi ) - m( \xi_1 ) | \leq C N^{-s} {{| \xi_1 |}^s} 
{{| \xi_1 |}^{- \theta}} {{| \xi_2 |}^\theta}
\end{equation}
for $0 \leq \theta \leq 1$. Recall that $m$ was defined to be smooth and 
monotone in \eqref{Idefined}.

Therefore, upon defining ${\mathcal{F}}(\nabla^\theta f)( \xi )
= {{|\xi|}^\theta} {\widehat{f}} ( \xi )$, we can write
\begin{equation*}
| {\mathcal{F}}
( \partial_x \{ I( u_h v_{vl} ) - I( u_h ) v_{vl} \}) ( \xi ) |
\leq |{\mathcal{F}}(\partial_x (\nabla^{-\theta} I(u_h )  
( \nabla^\theta {v_{vl}} ) ) ( \xi ) |.
\end{equation*}

We now estimate the left side of the bilinear estimate in this interaction
by
\begin{equation}
\X 0 {\half+} { \partial_x (\nabla^{-\theta} I (u_h ) )
( \nabla^\theta v_{vl} ) } 
\end{equation}
and by the bilinear estimate of Kenig, Ponce and Vega
\begin{equation}
\leq C \X 0 {\half+ } {\nabla^{-\theta} I (u_h) }  
\X 0 {\half+ } {\nabla^\theta v_{vl} } .
\end{equation}
The frequency support of $v_{vl}$ shows that $\X 0 {\half+ } {\nabla^\theta
v_{vl} } \lesssim \X 0 {\half+ } {v_{vl}} $. A moments thought shows
\begin{equation}
\X 0 {\half+ } {\nabla^{-\theta} I(u_h) } \leq N^{-\theta} 
\X 0 {\half+ } {I(u_h)}
\end{equation}
and the claim of the Proposition follows for the (very low)(high) interaction
by choosing $\theta > \frac{3}{4}$.

\noindent{\bf{Low/high interaction}}

The preceding calculations reduce matters to controlling
\begin{equation}
\X 0 {\half+ } {\partial_x \nabla^{-\theta} I(u_h) \nabla^\theta v_l }
\label{reduced}
\end{equation}
and we know that ${\widehat{u_h}}$ and ${\widehat{v_l}}$ are supported outside
the very low frequencies.

\begin{lemma}
Assume ${\widehat{u}} $ and ${\widehat{v}}$ are supported outside 
$\{ |\xi| < 1 \}$. Then
\begin{equation}
\X \alpha {-\half+ } {\partial_x (u v) } \leq C \X {-\gamma_1 } {\half+ } u
\X {-\gamma_2 } {\half+} v
\label{bilinear}
\end{equation}
provided
\begin{eqnarray*}
\alpha - (\gamma_1 + \gamma_2 ) &<& \frac{3}{4}, \\
\alpha - \gamma_i &<& \half,~i=1,2.
\end{eqnarray*}
\end{lemma}

We will apply the lemma momentarily with $\alpha =0, \gamma_1 = \gamma_2
= - \frac{3}{8}+$.

The proof of the lemma is contained in the proof of Theorem 2 in
\cite{CST99}. In particular, the support properties on ${\widehat{u}}, ~
{\widehat{v}}$ reduce matters to considering Cases A.3, A.4, A.6, B.3,
B.4, B.5 and B.6 in \cite{CST99}. The restriction $\alpha - (\gamma_1 +
\gamma_2 ) < \frac{3}{4}$ arises in Case A.4.c.ii of \cite{CST99} 
which is the region
containing the counterexample of \cite{KPVBilin}. Case B.4.b of \cite{CST99}
requires the other condition $\alpha - \gamma_i < \half.$

The lemma applied to \eqref{reduced} gives
\begin{equation*}
\leq C \X {-\frac{3}{8} +} {\half+} {\nabla^{-\theta} I(u_h ) }
\X {-\frac{3}{8}+ } {\half +} {\nabla^\theta v_l } .
\end{equation*}
Setting $\theta = \frac{3}{8} -$ leaves
\begin{equation*}
C \X 0 {\half+} {\nabla^{- \frac{3}{4}+} I(u_h ) } \X 0 {\half+ } {v_l }
\leq C N^{- \frac{3}{4}+}  \X 0 {\half+ } {I(u_h )} \X 0 {\half+} {v_l} 
\end{equation*}
which was to be shown.

\noindent{\bf{High/high interaction}}

In this region of the interaction, we do not take advantage of any cancellation
and estimate the difference with the triangle inequality 
\begin{equation*}
\X 0 {-\half + } {\partial_x \{ I(u_h ) I( v_h ) \} } +
\X 0 {-\half +} {\partial_x \{ I( u_h v_h ) \} }.
\end{equation*}
For the first contribution we use the lemma to get
\begin{equation}
\X {- \frac{3}{8}+} {\half+} {I(u_h) }
\X {-\frac{3}{8}+} {\half+} {I(v_h) } \leq N^{-\frac{3}{4}+} 
\X 0 {\half+} {I( u_h )}
\X 0 {\half+} {I( v_h ) }.
\end{equation}
The second contribution is bounded by throwing away $I$ and applying the
lemma,
\begin{eqnarray*}
\X 0 {-\half +} { \partial_x \{ u_h v_h \} } &\leq&
\X {- \frac{3}{8} +} {\half + } {u_h }
\X {- \frac{3}{8} +} {\half + } {u_h } \\
&\leq & N^{-\frac{3}{8} + s +} \X s {\half + } {u_h }
N^{-\frac{3}{8} + s +} \X s {\half + } {v_h } \\
&\leq & N^{-\frac{3}{4}+} \X 0 {\half +} {u_h}
\X 0 {\half +} {v_h } .
\end{eqnarray*}

% \bibliographystyle{plain}
% \bibliography{/home/u1/vis/colliand/Papers/Biblio/master}
%\bibliography{/home/colliand/Math/Papers/Biblio/master}

\begin{thebibliography}{10}

\bibitem{B1}
J.~Bourgain.
\newblock {F}ourier transform restriction phenomena for certain lattice subsets
  and applications to nonlinear evolution equations {I,II}.
\newblock {\em Geom. Funct. Anal.}, 3:107--156, 209--262, 1993.

\bibitem{Brefine}
J.~Bourgain.
\newblock {R}efinements of {S}trichartz' inequality and applications to
  2{D}-{N}{L}{S} with critical nonlinearity.
\newblock {\em International Mathematical Research Notices}, 5:253--283, 1998.

\bibitem{BAMS}
J.~Bourgain.
\newblock {\em Global solutions of nonlinear {S}chr\"odinger equations}.
\newblock American Mathematical Society, Providence, RI, 1999.

\bibitem{CKSTTDNLS}
J.~Colliander, M.~Keel, G.~Staffilani, H.~Takaoka, and T.~Tao.
\newblock {G}lobal well-posedness for {S}chr\"odinger equations with derivative.
\newblock (preprint), 2001.

\bibitem{CKSTTNLS}
J.~Colliander, M.~Keel, G.~Staffilani, H.~Takaoka, and T.~Tao.
\newblock {G}lobal well-posedness of 2d {N}{L}{S}.
\newblock (in preparation), 2001.

\bibitem{CKSTTKdV}
J.~Colliander, M.~Keel, G.~Staffilani, H.~Takaoka, and T.~Tao.
\newblock {S}harp {G}lobal well-posedness of periodic and nonperiodic generalized
  {K}orteweg-de {V}ries equations.
\newblock (in preparation), 2001.

\bibitem{CST99}
J.~E. Colliander, G.~Staffilani, and H.~Takaoka.
\newblock {G}lobal wellposedness of {K}d{V} below $L^2$.
\newblock {\em Mathematical Research Letters}, 6(5,6):755--778, 1999.

\bibitem{KeeTao98}
M.~Keel and T.~Tao.
\newblock {L}ocal and {G}lobal {W}ell-{P}osedness of {W}ave {M}aps on
  ${\R}^{1+1}$ for {R}ough {D}ata.
\newblock {\em International Mathematical Research Notices}, 21:1117--1156,
  1998.

\bibitem{KPVBilin}
C.~Kenig, G.~Ponce, and L.~Vega.
\newblock A bilinear estimate with applications to the {K}d{V} equation.
\newblock {\em J. Amer. Math. Soc.}, 9:573--603, 1996.

\bibitem{Staff97}
G.~Staffilani.
\newblock On the growth of high {S}obolev norms of solutions for {K}d{V} and
  {S}chr\"odinger equations.
\newblock {\em Duke Math. J.}, 86(1):109--142, 1997.

\end{thebibliography}

\enddocument